\documentclass[12pt]{article}
\usepackage{graphicx}
\usepackage{amsmath,comment,fancyhdr,color,graphics}
\usepackage{amsmath}
\usepackage{amsfonts}

\def\N{\mathbb{N}}

\def\0{{\bf 0}}

\def\N{\mathbb{N}}

\def\0{{\bf 0}}

\def\real{\mathbb{R}}

\newcommand{\divc}{\mbox{dim}_{\mbox{\tiny VC}}}
\newcommand{\devc}{\mbox{dens}_{\mbox{\tiny VC}}}
\newcommand{\didu}{\mbox{dim}^*}
\newcommand{\dedu}{\mbox{dens}^*}
\newtheorem{theo}{Theorem}
\newtheorem{coro}{Corollary}
\newtheorem{definition}{Definition}

\newtheorem{example}{Example}
\newtheorem{fact}{Fact}
\newtheorem{claim}{Claim}
\newtheorem{lemma}{Lemma}

\newcommand{\eps}{\epsilon}
\newcommand{\beq}{\begin{equation}}
\newcommand{\enq}{\end{equation}}
\newcommand{\beqa}{\begin{eqnarray}}
\newcommand{\enqa}{\end{eqnarray}}

\newcommand{\neno}{\newline\noindent}

\def\N{\mathbb{N}}

\newcommand{\A}{{\cal A}}
\renewcommand{\H}{{\cal H}}

\newcommand{\M}{{\cal M}}

\newcommand{\F}{{\cal F}}
\newcommand{\G}{{\cal G}}

\newcommand{\T}{{\cal T}}

\newcommand{\dd}{\mbox{d}}
\def\qed{\hfill\hbox{${\vcenter{\vbox{
        \hrule height 0.4pt\hbox{\vrule width 0.4pt height 6pt
        \kern5pt\vrule width 0.4pt}\hrule height 0.4pt}}}$}}

\begin{document}

\begin{center}
\begin{Large}
Metric Entropy estimation using o-minimality Theory
\end{Large}
\vskip 2pc
\begin{large}
A. Onshuus and  A. J. Quiroz\footnote{Dpto. de
Matem\'aticas, Universidad de Los Andes. Address: Dpto. de Matem\'aticas, Universidad de Los Andes, Carrera 1, Nro. 18A-10, edificio H, Bogot\'a, Colombia. Phone: (571)3394949, ext. 2710. Fax: (571)3324427. 
e-mails: aonshuus@uniandes.edu.co, aj.quiroz1079@uniandes.edu.co}
\end{large}
\end{center}
\vskip 1pc
\begin{small}
{\bf Abstract}
It is shown how tools from the area of Model Theory, specifically from the Theory of o-minimality, can be used to prove that a class of functions is VC-subgraph (in the sense of \cite{dudc}), and therefore satisfies a uniform polynomial metric entropy bound. We give examples where the use of these methods significantly improves the existing metric entropy bounds.  The methods proposed here can be applied to finite dimensional parametric families of functions without the need for the parameters to live in a compact set, as is sometimes required in theorems that produce similar entropy bounds (for instance Theorem 19.7 of \cite{vdv}).
\end{small}
\section{Introduction}
VC-dimension and metric entropy are fundamental concepts in modern asymptotic statistics and the theory of statistical learning, due to their applicability in establishing uniform convergence results, such as Uniform Laws of Large Numbers and Uniform Central Limit Theorems (see \cite{dudc}, \cite{poll}, \cite{vdv} or \cite{vdvw}). The metric entropy of a class of functions is a measure of the size of the class, and, for classes of functions, it plays a roll, with regards to asymptotics, very similar to the one played by VC-density (or VC-dimension) for collections of sets. Finding a tight upper bound for the VC-density (or VC-dimension) of a class of sets
or for the metric entropy of a class of functions (when this is small enough) will be useful in establishing speed of convergence in the Uniform Law of Large Numbers, as shown in Section 2.6 of \cite{poll} or in the Central Limit Theorem for certain functionals of the empirical process, as explained in Section 3.4 of \cite{vdvw}.

In the present article, we will show how working with the VC-density of an associated class of sets and taking advantage of recent results from the theory of o-minimality, allows for significantly improving bounds on the metric entropy of certain classes of functions and for finding tight entropy bounds for certain classes for which other methods would not work.

Even though most of the concepts we have mentioned are by now classical, we briefly review, in the following subsection, their definitions and their use in theorems of asymptotic statistics. This will be useful for setting our notation and establishing a context for the calculations in the following section.

\subsection{VC theory in asymptotic statistics}
\begin{definition}
Given a collection of measurable sets in $d$-dimensional Euclidean space, $\A$, and a finite set $F\subset\real^d$, let
\[
F\cap\A=\{F\cap A: A\in \A\} \> \mbox{ and }\> \Delta_{\A}(n)=\sup_{\{F:|F|=n\}}|F\cap\A|,
\]
where $|\cdot|$ denotes cardinality of a (usually finite) set.\footnote{In other places $|\cdot|$ will denote absolute value, as usual, but (we hope) this will cause no confusion since the meaning should be clear from the context in each case.} When $\Delta_{\A}(n)$ is bounded by a polynomial in $n$, the class $\A$ is said to be a Vapnik-Cervonenkis class or, shortly, a VC-class. The VC-density of $\A$, $\devc(\A)$, is the infimum of the set of positive reals, $r$, such that a constant $C>0$ exists (possibly depending on $r$) such that $\Delta_{\A}(n)\leq C\,n^r$ for all $n\in \N$. The VC-dimension of $\A$, $\divc(\A)$, is the largest positive integer $m$ such that $\Delta_{\A}(m)=2^m$. If no such $m$ exists,  $\divc(\A)=\infty$.
\end{definition}

Under measurability conditions, VC-classes satisfy a non-parametric Uniform Law of Large Numbers, in the sense that if $\A$ is a VC-class, for i.i.d. data $X_1,\dots,X_n$, sampled from a probability distribution $P$ on $\real^d$ and with
\beq\label{r1}
P_n(A)=\frac{|\{i:1\leq i\leq n,\>X_i\in A\}|}{n},
\enq
for a set $A$, then
\[
\sup_{A\in\A}|P_n(A)-P(A)|\to 0,\>\>\mbox{a.s., as } n\to\infty.
\]
Examining the proof of this Uniform Strong Law reveals that the value of the VC-density is what is actually involved in the arguments leading to this result. Still, most often, asymptotic statisticians have resorted to bounds on the VC-dimension to establish this type of strong laws. This is probably due, at least in part, to the following facts:
\neno (i) Finiteness of the VC-dimension is equivalent to finiteness of the VC-density and
\neno (ii) During the 1980's, several methods were developed for bounding the VC-dimension of a class of sets.

A notion related to VC-dimension and density, that will be used below, is the dual dimension of Assouad, \cite{assd}, defined as follows.
\begin{definition}
For a class of sets in $d$-dimensional Euclidean space, $\A$, and a finite sub-collection, $\H\subset \A$, let $\mbox{At}(\H)$ denote the set of atoms of the finite algebra generated by $\H$. For $n\in\mathbb{N}^+$, let
$\Delta_{\A}^*(n)=\sup\,\{|\mbox{At}(\H)|: \H\subset\A,\,|\H|=n\}$. Assouad's dual dimension is
\beq\label{r3}
\didu(\A)=\sup\,\{m\in\mathbb{N^+}:\mbox{ there exists }\H\subset\A, \mbox{ with }|\H|=m,\,|\mbox{At}(\H)|=2^m\}.
\enq
The dual density of $\A$, $\dedu(\A)$ is the infimum of the positive reals, $r$, such that for some constant $C>0$, $\Delta_{\A}^*(m)\leq C\,m^r$ for every $m\in\mathbb{N}^+$.
\end{definition}

When it comes to asymptotic results over classes of functions, the concept of metric entropy plays a role similar to that of VC-density for classes of sets.

\begin{definition} Let $\F\subset L^p(Q)$ for $p=1$ or 2, and a probability measure $Q$ on $\real^d$.
For $\eps>0$, the $\eps$-covering number of $\F$ with respect to $Q$, $N_p(\eps,\F,Q)$, is the minimum natural $m$ such that
there exist functions $g_1,g_2,\dots,g_m\in L^p(Q)$ satisfying that, for every $f\in \F$, there is a $j\in \{1,\dots,m\}$ such that $\|f-g_j\|_{p,Q}<\eps$ where $\|\cdot\|_{p,Q}$ is the norm of $L^p(Q)$.
$H_p(\eps,\F,Q)=\log N_p(\eps,\F,Q)$ is called the metric entropy of $\F$.
\end{definition}

In order to state a law of large numbers over $\F$, let again $X_1,\dots,X_n$ denote an i.i.d. sample from a probability distribution $P$ on $\real^d$, and, for each integrable function $f$, let
$P_n\,f=(1/n)\sum_{i\leq n}f(X_i)$, be the empirical integral of $f$, while $Pf=\int f(x)\dd P(x)$.  The class $\F$ is said to have an envelope function $F\in L^p(Q)$ whenever $|f(x)|\leq F(x)$ for all $f\in \F$ and every $x\in\real^d$. A Uniform Law of Large Numbers holds over $\F$ with respect to $P$, when
\beq\label{r4}
\sup_{f\in\F}|P_n(f)-Pf|\to 0,\mbox{ a.s. as } n\to\infty.
\enq
Different results exist in the literature connecting bounds on the metric entropy of a class to Uniform Laws as (\ref{r4}).

Suppose that the $L^p$ covering number, $p=1,2$, of the class $\F$ with envelope function $F$ is small enough as to satisfy a polynomial bound such as
\beq\label{r13}
\sup_QN_p(\eps\|F\|_{p,Q},\F,Q)\leq A\left(\frac{1}{\eps}\right)^{B}
\enq
where $A$ and $B$ are positive constants and the bound is uniform over all choices of the probability measure $Q$. In this case we will say that {\it the class $\F$ has polynomial $L^p$ covering number (with exponent $B$)}. When a class of functions has polynomial covering number, more things can be said regarding asymptotics. When the  $L^1$ covering number is polynomial and the envelope function $F$ is bounded, then the uniform strong law (\ref{r4}) can be improved with a uniform speed of convergence:
\beq\label{r14}
\sup_{f\in\F}|P_n(f)-Pf|\ll \frac{\log n}{\sqrt{n}},\mbox{ a.s. as } n\to\infty.
\enq
where $a_n\ll b_n$ means that $a_n/b_n\to 0$, as can be deduced from Theorem 37 in \cite{poll}. Similarly, for classes
with polynomial $L^2$ covering number, results on the speed of convergence in the Central Limit Theorem for the Empirical Process over a class of functions can be obtained, as explained in Section 3.4 of \cite{vdvw}. Thus, establishing the polynomial $L^p$ covering number property for a class of functions is quite relevant from the asymptotic viewpoint.

In the current literature, there exist two ways of proving that a class of functions $\F$ has polynomial $L^p$ covering number. One is through the notion of VC-subgraph classes, to be discussed in a moment. The other is through the total boundedness of the finite dimensional set of parameters that index the functions in $\F$. This second method appears as Example 19.7 in \cite{vdv} and we will refer to it in the sequel as the {\it bounded parameter space} method. The first method has the advantage of not needing the parameter space to be totally bounded and other technical (smoothness) conditions required in the bounded parameter space method.
The purpose of the main result in the present article is to significantly simplify the verification that a parametric class of functions is VC-subgraph, thus obtaining that the class has polynomial $L^p$ covering numbers, for $p=1$ and 2. In examples we will show several classes that appear in concrete applications to be VC-subgraph.

First, let us recall the definition of VC-subgraph classes introduced in \cite{poll82}, although this name comes from \cite{dudc}.
\begin{definition}
For a class of functions on $\real^d$, $\F$, and $f\in\F$, the subgraph of $f$ is the set
\beq\label{r7}
\mbox{subgraph}(f)=\{(x,t)\in\real^{d+1}:0\leq t\leq f(x)\mbox{ or }0>t>f(x)\}.
\enq
The class of all subgraphs of functions in $\F$, \mbox{subgraph}$(\F)$, is a collection of sets in $\real^{d+1}$. When \mbox{subgraph}$(\F)$ is a VC-class, $\F$ is called a VC-subgraph class.
\end{definition}
Careful reading of the proof of Lemma 25 in Chapter 2 of \cite{poll} gives the following:
\begin{theo}\label{VC versus entropy}
If $\F$ is a VC-subgraph class with envelope $F\in L^p(Q)$ and \newline\noindent
$r=\devc(\mbox{subgraph}(\F)$), then, for any $\eta>0$
\beq\label{r5}
N_p(\eps\|F\|_{p,Q},\F,Q)\leq A\left(\frac{1}{\eps}\right)^{r+\eta}
\enq
where the constant $A$ depends only on $r$ and $\eta$ (not on $Q$). That is, $\F$ has polynomial $L^p$ covering number with exponent $r+\eta$, for every positive $\eta$.
\end{theo}

The next subsection includes some facts found recently in the context of Model Theory to be used later.

\subsection{Some definitions and results from o-minimality}

O-minimality and results we will be using are a subarea of model theory, in the sense of mathematical logic. We will need some definitions,
although we will try to give enough examples so that the un-familiarized reader can get an idea of the concepts we will need.

\subsubsection{First order logic}

In order to work in model theory, one fixes a language $\mathcal L$ (say, the language of rings with unity $\mathcal{ L}:=\{+, \cdot, 0, 1\}$)
and a structure which interprets each symbol that appears in $\mathcal L$ (for example the real field and the complex fields are both structures in
the language of rings with unity) which will be called an $\mathcal L$-structure.

In this paper we will always work with a language $\mathcal L$ which includes the language of ordered rings $\mathcal{ L}_{o. ring}:= \left\{+, \cdot, 0, 1, <\right\}$, and
with structures whose universe is the real numbers as an ordered field, fixing an interpretation of the symbols in $\mathcal L$, which will
usually be the standard interpretation. For example, the structure with universe $\mathbb R$, associated with the language $\mathcal L:=\left\{+, \cdot, 0, 1, <, e^x\right\}$, will be the real numbers, with the natural interpretation of, respectively, the addition, multiplication, additive identity, multiplicative identity, order, and exponential function.

\begin{definition}
Let $\mathcal L$ be a fixed language which includes $\mathcal L_{o. ring}$, and let $(\mathbb R, +, \cdot, 0, 1, <, \dots)$ be an
$\mathcal L$-structure (here the dots stand for whichever relations or function symbols we want to add to the language).
An \emph{$\mathcal L$-definable subset of $\mathbb R$} is the set of realizations in our structure $(\mathbb R, +, \cdot, 0, 1, <, \dots)$
of a formula which uses only symbols from $\mathcal L_{o. ring}$, besides the logic symbols $=, \vee, \wedge, \Rightarrow, \Leftrightarrow, \neg \text{(the symbol for ``it is not true that'')},
\forall $ and $\exists$. \end{definition}

Given an $\mathcal L$-structure $\mathcal M$, an $\mathcal L$-formula $\phi(\bar x, \bar y)$ and $\mathcal M$-tuples $\bar a$ and $\bar b$, we will say that
\[
\mathcal M\models \phi(\bar a, \bar b)
\]
if the formula $\phi(\bar a, \bar b)$ is true in $\mathcal M$.

\begin{example}
Let $\phi(x,y)$ be the formula $\exists z\  z^2=(y-x)$, and let $\mathbb Q$ and $\mathbb R$ be the rational and the real fields, respectively. Then
\[
\mathbb R\models \phi(1, 3)
\]
but
\[
\mathbb Q\not\models \phi(1, 3).
\]
\end{example}

Of course, if a formula $\phi$ has no free variables (so that all variables appearing in $\phi$ are quantified by either $\exists$ or $\forall$) we don't need to
replace any variables to know the truth or falsehood of $\phi$ in any structure of the language. Such formulas are called \emph{sentences}. So for
example, $\phi:=\forall x \exists y \ y\cdot y=x$ is a sentence true in the complex field, false in the real field, but true in the structure
$(\mathbb R^{\geq 0}, +, \cdot, 0, 1)$.

\begin{example}
$\left.\right.$
\begin{itemize}
\item The unit disk in $\mathbb R^2$   is definable in the structure
$(\mathbb R, +, \cdot, 0, 1, <)$ by the formula $x^2+y^2<1$.

\item The set of integers are a definable subset of $\mathbb R$ in the structure $(\mathbb R, +, \cdot, 0, 1, <, sin(x))$, since it is the set or realizations of the formula $\exists y, \ sin(y)=0 \wedge x\cdot \pi=y$.

\item By Fact \ref{R o-min} below, the integers are not definable in the structure $(\mathbb R, +, \cdot, 0, 1, <)$.

\item The derivative of a function $f(x)$ is definable in the structure $(\mathbb R, +, \cdot, 0, 1, <, f(x))$ by replacing $f'(x)=y$ by the formula

\[
\phi(x,y):= \forall \epsilon   \exists \delta
\forall h \left( -\delta<h<\delta\right)\Rightarrow \left(\left| f\left(x+h\right)-f\left(x\right)-h y \right|< \left|h\epsilon\right| \right).
\]
(Here the absolute value can be defined in the standard way or, since by $|x|:=\sqrt {x^2}$ we can define $|x|=y$ by the formula $\theta(x,y):=\left(y^2=x^2\right) \wedge \left(y\geq 0\right)$.)

\end{itemize}
\end{example}

\subsubsection{o-minimality}
The main logic definition of this paper is the following.

\begin{definition}
Let $\mathcal L$ be a fixed language which includes $\mathcal L_{o. ring}$, and let $(\mathbb R, +, \cdot, 0, 1, <, \dots)$ be an
$\mathcal L$-structure. We will say that $(\mathbb R, +, \cdot, 0, 1, <, \dots)$ is o-minimal if and only if every $\mathcal L$-definable subset of
$\mathbb R$ is a finite union of open intervals and points.
\end{definition}

The following is a well known theorem of Tarski (see \cite{Ta}) (not stated originally in this precise manner, since the concept of o-minimality came later).

\begin{fact}\label{R o-min}
$(\mathbb R, +, \cdot, 0, 1, <)$ is o-minimal.
\end{fact}

Notice that the subset $\mathbb Z$ of $\mathbb R$ is not a finite union of intervals and points, so in particular o-minimality implies
it is not definable in the real field (or in any o-minimal expansion of the real field).

\medskip

The study of o-minimal theories started without any other examples of o-minimal expansions of the real field. But it really
started becoming a main area of model theory with the following theorem due to Wilkie (\cite{Wi}).
\begin{fact}\label{exp o-min}
$\mathbb R_{exp}:=(\mathbb R, +, \cdot, 0, 1, <, e^x)$ is o-minimal.
\end{fact}

This was later generalized by van den Dries, Macintyre and Marker (\cite{vdDMM}) to the following statement:

\begin{fact}\label{exp o-min an}
Let $R_{exp, an}$ be the real field expanded by the exponential functions and a function symbol for every
analytic function with domain $[-1, 1]^m$ for some $m$. Then, $R_{exp, an}$ is o-minimal.
\end{fact}
(For example, even though $(\mathbb R, +, \cdot, 0, 1, <, sin(x))$ is not o-minimal, the structure
\[\left(\mathbb R, +, \cdot, 0, 1, <, sin\left(x\right)|_{[-1,1]}\right),\] where $sin(x)|_{[-1,1]}$ is the restriction of
the sine function to the closed interval $[-1,1]$, is o-minimal.

This, together with a result of Speissegger, will cover most of the examples we will consider in this paper. But in order to state Speissegger's result we need
the following definition.

\begin{definition}
Let $(\mathbb R, +, \cdot, 0, 1, <, \dots)$ be any expansion of the real field. We will say that a differential
equation is \emph{Pfaffian over $(\mathbb R, +, \cdot, 0, 1, <, \dots)$} if it is given by a system of equations
of the form
\[
\frac{\partial f_i}{\partial x_j} =P_{i,j}\left(\bar x, f_1\left(\bar x\right), \dots, f_{i}\left(\bar x\right)\right)
\]
where $f_i\left(\bar x\right)$ and $P_{i,j}(\bar y)$ are definable functions in $(\mathbb R, +, \cdot, 0, 1, <, \dots)$, $j$ varies through the number of variables and $1\leq i\leq N$ for some positive integer $N$.
\end{definition}

The following is due to Speissegger (\cite{Sp}):
\begin{fact}\label{Speissegger}
Let $(\mathbb R, +, \cdot, 0, 1, <, \dots)$ be any o-minimal expansion of the real field, and let $f(\bar x)$ be
the solution of a Pfaffian differential equation in $(\mathbb R, +, \cdot, 0, 1, <, \dots)$. Then the structure
$(\mathbb R, +, \cdot, 0, 1, <, \dots, f)$ is o-minimal.
\end{fact}

Notice that since $e^x$ is a solution of $\frac{\partial f}{\partial x}=f$, Wilkie's result follows from the o-minimality of the real field
$(\mathbb R, +, \cdot, 0, 1, <)$ and Speissegger's result.

It follows that the Pfaffian closure of $\mathbb R_{exp, an}$ (which we will denote $\mathbb R_{an, Pfaff}$) is o-minimal\footnote{Here
we use closure in a manner analogous to ``algebraic closure'': a structure $\mathcal M$ is Pfaffian closed if given any function $f$, if $f$ is Pfaffian over $\mathcal M$, then
$f$ is definable in $\mathcal M$. The Pfaffian closure of $\mathcal R$ is the smallest structure containing $\mathcal R$ which is Pfaffian closed}.

\subsubsection{Uniform definable families of sets}

We will begin with a definition. By a \emph{Uniform definable family of definable sets in $\mathcal M$}
we mean a family of definable subsets, all of which are given by changing the parameters in a fixed formula in the language $\mathcal L$.
Formally,

\begin{definition}
Let $\mathcal L$ be any language and let $\mathcal M:=(M,\dots)$ be an $\mathcal L$-structure. We will say that $\mathcal F$ is a
uniform definable family of definable subsets of $M^n$ if there is an $\mathcal L$-formula $\phi(\bar x; \bar y)$ such that
\[\mathcal{ F}:=\left\{X_{\bar b}\right\}_{\bar b\in M^d}\] where
\[X_{\bar b}:=\left\{ \bar a\in M^n : \ \phi(\bar a; \bar b)\text{ is true in $\mathcal M$}\ \right\}.\]
The tuples $\bar b$ vary in $M^d$ (where $d$ is the dimension of the variable $\bar y$ in the formula $\phi(\bar x; \bar y)$) and will be called
the ``parameters'' of the subset $X_{\bar b}$.
\end{definition}

\begin{example}\label{Example 2}
$\left.\right.$
\begin{itemize}
\item Since the semi-spaces of $\mathbb R^n$ are all definable by a formula
\[b_1\cdot x_1+b_2\cdot x_2+\dots+b_n\cdot x_n+b_{n+1}<0,\]
they are a uniform definable family in $(\mathbb R, +, \cdot, 0, 1, <)$.

\item The family $\mathcal A:=\{X_\lambda\}_{\lambda \in \mathbb R}$ where \[X_\lambda:=\left\{x : x\geq 0, \ 0\leq \frac{x^\lambda-1}{\lambda}\right\}\] is
uniformly definable in the structure
$(\mathbb R, +, \cdot, 0, 1, <, e^x)$.
(Technically, we would need to replace $x^\lambda-1=z$ with the formula $\exists y\ \left( e^{\lambda \cdot y}-1=z\right)\ \wedge\  \left(e^y=x\right)$, but this is all definable by a first order formula.)
\end{itemize}
\end{example}

\subsubsection{The main theorem}

The main theorem relating o-minimality to VC-density, was explicitly stated and proved in \cite{ADHMS}, although the result is already contained in the paper
\cite{KMac}. Here we state the result, followed by an immediate implication in terms of the VC-density of the subgraphs of a class of functions.

\begin{theo}\label{ADHMS}
Let $\mathcal R:=(\mathbb R, +, \cdot, 0, 1, <, \dots)$ be an o-minimal expansion of the real field, and let $\mathcal F:=\{X_{\bar b}\}_{\bar b\in \mathbb R^d}$  be a uniform definable
family of sets defined by the formula $\phi(\bar x; \bar y)$ with $\bar x$ an $m$-tuple of variables and $\bar y$ a $d$-tuple of variables.
\[
\mathcal F:=\{X_{\bar b}\}_{\bar b\in \mathbb R^d}:=\left\{ \left\{\bar a\in \mathbb R^m :  \mathbb R\models \phi\left(\bar a, \bar b\right)\right\}:\bar b\in \mathbb R^d\right\}.
\]
Then the VC-density of $\mathcal F$ is at most $d$.
\end{theo}

It follows for instance that, since the family $\mathcal
A:=\{X_{\lambda}\}_{\lambda \in \mathbb R}$ in Example \ref{Example 2} is a one parameter family
uniformly defined in the o-minimal structure $\mathbb R_{exp}$, the
VC-density of $\mathcal A$ is at most 1. Recall that by definition this means
there is some real constant $C>0$ such that $\Delta_{\mathcal
A}(n)<C\cdot n$ for all $n$.

The proof of Theorem \ref{ADHMS} in
\cite{ADHMS} is done by induction on the length of the parameter
set, and the proof for the 1-case will actually give us a bound for
$C$ in this particular example (which we will work out in Subsection
\ref{2.2}).

More generally, because the subgraph of the function $f(\bar x)$ is the set
\[\left\{\left(\bar x, y\right) : 0\leq y\leq f(\bar x)\right\}\cup \left\{\left(\bar x, y\right) : 0\geq y\geq f(\bar x)\right\}\] we can state the following
general result about the VC-density of subgraphs of uniformly definable functions:

\begin{coro}
Let $\mathcal R:=(\mathbb R, +, \cdot, 0, 1, <, \dots)$ be an o-minimal expansion of the real field, and let $\mathcal F:= \{f_{\bar b}(\bar x)\}_{\bar b\in \mathbb R^d}$  be a uniform definable
family of functions defined by the formula $\phi(y, \bar x; \bar z)$ with $\bar x$ an $m$-tuple of variables and $\bar z$ a $d$-tuple of variables. Explicitly,
$\mathcal F:=\{f_{\bar b}(\bar x)\}_{\bar b\in \mathbb R^d}$ is defined so that for any $\bar x$ and $y$, we have $f_{\bar b}(\bar x)=y$ if and only if $\phi(y, \bar x; \bar b)$ holds.

Then the VC-density of $\mbox{subgraph}({\mathcal F})$ is at most $d$.
\end{coro}

The following is a direct consequence of Theorems \ref{VC versus entropy} and \ref{ADHMS} and it is the tool proposed here for statistical applications.

\begin{coro}\label{Metric}
Let
\[
\mathcal F:=\{X_{\bar b}\}_{\bar b\in \mathbb R^d}:=\left\{ \left\{\bar a\in \mathbb R^m :  \mathbb R\models \phi\left(\bar a, \bar b\right)\right\} : \bar b\in \mathbb R^d \right\}.
\]
be a (parametric) family of functions on $\mathbb R^m$ (uniformly) definable in an o-minimal structure with $d$ parameters, and assume also that $\mathcal F$ has bounded envelope function $F$. Then,  $\F$ has polynomial $L^p$ covering number with exponent $d+\eta$, for any $\eta>0$ and $p=1, 2$.
\end{coro}

\section{Bounding the metric entropy of certain classes of functions}\label{2.2}

Next, we consider certain classes of functions that have appeared in the statistical literature and show how to improve the bounds that have been reported on their metric entropy.

\subsection{Transformations to elliptical symmetry}
Our first example appeared in \cite{npq} in connection to the estimation of transformations of multivariate data to elliptical symmetry. In order to establish the efficiency of the method proposed there, part of the problem reduces to the consideration of the class of functions on $\real^+$, $\T$, defined by
\[
T_{\lambda}(x)=\frac{x^{\lambda}-1}{\lambda}, \mbox{ for }x\in\real^+   \qquad\mbox{ and }\qquad  \T=\{T_\lambda:\lambda \in \Lambda\},
\]
where $\Lambda$ is a bounded interval. In order to study the class subgraph($\T$), its dual class, subgraph($\T)^*$, formed by the sets
\[
T^{dual}(x,t)=\{\lambda\in\Lambda:\> 0\leq t\leq\frac{x^{\lambda}-1}{\lambda}\>\>\mbox{ or }\>\>
\frac{x^{\lambda}-1}{\lambda}\leq t<0\}
\]
for $(x,t)\in (\real^+\times\real)$, was considered. Since (as we shall see in the proof of Lemma 1) each $T^{dual}(x,t)$ is the union of at most two intervals, it was argued in \cite{npq} that the VC-dimension of subgraph($\T)^*$ is bounded by 4.\footnote{Taking points 1,2,3,4 and 5, one can not, with two intervals, pick the set \{1,3,5\}.} Then, by Proposition 2.13 in
\cite{assd}, it follows that the VC-dimension of subgraph($\T$) is bounded by $2^4=16$ and, therefore, $\devc(\mbox{subgraph}(\T))$ will be bounded by 16.

Since the subgraphs of $T_\lambda$ are a uniformly definable family in the o-minimal structure $\mathbb R_{exp}$,
by Theorem \ref{ADHMS} we know that the VC-density of the subgraphs of the functions $T_{\lambda}$ is bounded by the size of the parameter set, so it is bounded by a linear function. Now, a closer analysis of the methods in the one dimensional case of the proof of Theorem \ref{ADHMS} in \cite{ADHMS}, will give a precise bound for this family which might be useful for getting precise bounds in any one dimensional set.

\begin{lemma}
For any fixed $\lambda\in \mathbb R$, consider the subgraph
\[
S_{\lambda}=\left\{(x,t)\in \mathbb R^{\geq 0}\times \mathbb R : 0\leq t\leq \frac{x^{\lambda}-1}{\lambda}\vee 0\geq t\geq \frac{x^{\lambda}-1}{\lambda}\right\}
\]
and let $\mathcal A:=\{S_\lambda\}_{\lambda\in \mathbb R}$. Then $\Delta_{\mathcal A}(n)\leq n+1$.
\end{lemma}

\noindent \emph{Proof:}
We will again work with the dual subsets but in a different manner. Let
\[X_n:=\left\{\left(x_1, t_1\right), \left(x_2, t_2\right), \dots, \left(x_n, t_n\right)\right\}\] be any $n$ points in $\mathbb R^2$, and
we want to bound the number of sets in $X_n\cap \mathcal A$. Now, for each pair $(x_i, t_i)$ let
\[
T^{dual}(x_i, t_i):=\{\lambda : (x_i, t_i)\in S_{\lambda}\}.
\]

The next observation is trivial, but it is the central piece of our argument.
\begin{claim}
If $S_{\lambda_1}$ and $S_{\lambda_2}$ define different subsets of $X_n$, then for some $(x_i, t_i)$ we have that
\[
\lambda_1\in T^{dual}(x_i, t_i) \not\Leftrightarrow \lambda_2\in T^{dual}(x_i, t_i)
\]
so that $\lambda_1$ is in the set $T^{dual}(x_i, t_i)$ and $\lambda_2$ isn't, or viceversa.
\end{claim}

Notice that if
\[
\lambda\in \bigcap_{i\in I} T^{dual}\left(x_i, t_i\right)  \cap
\bigcap_{j\not\in I} \left(\mathbb R\setminus T^{dual}\left(x_j, t_j\right)\right),\] then
\[S_\lambda\cap X_n=\{(x_i,t_i)\}_{i\in I}.\]

It follows from the claim and the above observation (by an easy and insightful argument left to the reader)
that the number of sets  $X_n\cap \mathcal A$ is equal to the number of non empty intersections
of the sets $X_n\cap \mathcal A$ and their complements, so that
\[
\left| X_n\cap \mathcal A\right|=\left|\left\{ I\subset \left\{1, \dots n\right\} : \bigcap_{i\in I} T^{dual}\left(x_i, t_i\right)  \cap
\bigcap_{j\not\in I} \left(\mathbb R\setminus T^{dual}\left(x_j, t_j\right) \right) \neq \emptyset \right\} \right|.
\]

For notation purposes, for any $I\subset \left\{1, \dots n\right\}$, let
\[
B_I:=\bigcap_{i\in I} T^{dual}\left(x_i, t_i\right)  \cap
\bigcap_{j\not\in I} \left(\mathbb R\setminus T^{dual}\left(x_j, t_j\right) \right).
\]

So we need to count the subsets $I$ which give consistent (non-empty) boolean combinations.

For any $(x_i, t_i)$ let $c_i:= \inf (T^{dual}(x_i, t_i))$ if it exists ($-\infty$ otherwise)
and let $d_i:= \sup (T^{dual}(x_i, t_i))$ if it exists ($\infty$ otherwise). The derivative
of  \[
h(\lambda)=\frac{x^\lambda-1}{\lambda}
\]
as a function of $\lambda$, is always positive, as one can easily verify\footnote{In fact, the minimum of $h'(\lambda)$ is always 0, at $\lambda=0$.}. Furthermore,
for $x>1$, $h(\lambda)$ is always positive with infimum 0, while for $x<1$, $h(\lambda)$ is always negative with
supremum 0.
It follows that
$$
\begin{array}{ll}
 T^{dual}(x_i, t_i):=\left[c_i, \infty \right) & \text{ if } x_i>1, t_i\geq \ln(x_i), \\
 T^{dual}(x_i, t_i):=\left[c_i,0\right)\cup \left(0, \infty \right) & \text{ if } x_i>1, 0 < t_i<\ln(x_i), \\
 T^{dual}(x_i, t_i):=\mathbb{R}\setminus\{0\} & \text{ if } x_i>1, t_i= 0,\\
 T^{dual}(x_i, t_i):=\emptyset & \text{ if } x_i>1, t_i< 0,\\
 T^{dual}(x_i, t_i):=\emptyset & \text{ if } x_i<1, t_i>0, \\
 T^{dual}(x_i, t_i):=\mathbb{R}\setminus\{0\} & \text{ if } x_i<1, t_i= 0,\\
 T^{dual}(x_i, t_i):=\left(-\infty,0\right)\cup \left(0, d_i \right] & \text{ if } x_i<1, 0> t_i>\ln(x_i), \\
 T^{dual}(x_i, t_i):=\left(-\infty, d_i\right] & \text{ if } x_i<1, 0>\ln(x_i)\geq t_i,
\end{array}
$$
where all values of $c_i$ and $d_i$ are finite.
Let $I\subset \{1, 2, \dots, n\}$ be any subset for which $B_I\neq\emptyset$.

Working on the set $\mathbb{R}\setminus\{0\}$ (disregarding the zero element) and ignoring trivial values of the set $T^{dual}(x_i, t_i)$, we can assume that all our $T^{dual}(x_i, t_i)$ are of the form $[c_i,\infty)$ or $(-\infty,d_i]$. Assume, without loss of generality, that the $c_i$ are listed in increasing order and so are the $d_l$: For $i<i'$, $c_i<c_{i'}$ and for $l<l'$, $d_l<d_{l'}$. Write $T_i$ for $T^{dual}(x_i, t_i)$.

If $B_I$ is non empty, $i\in I$ and $T_i$ is of the form $[c_i,\infty)$, then for $i'<i$ and $T_{i'}=[c_{i'},\infty)$, we must have
$i'\in I$ (otherwise $B_I$ would be empty). Similarly, if $l\in I$ and $T_l$ is of the form $(-\infty,d_l]$, then for $l'>l$ and $T_{l'}=(-\infty,d_{l'}]$, we must have $l'\in I$.

In order to define a non empty $B_I$, we first choose $i_1$ as the largest $i$ such that $T_i=[c_i,\infty)$ and $i\in I$. Let $i_2$ be the next $i$ such that $T_i=[c_i,\infty)$.
Then $B_I$ must be contained in $[c_{i_1},c_{i_2})$. If $i_2$ does not exist, $B_I$ must be contained in $[c_{i_1},\infty)$. Let
$L(i_1)$ denote the set of indices $l$ such that $T_l$ is of the form $(-\infty,d_l]$ and $d_l\in [c_{i_1},c_{i_2})$. Then, $B_I$ is completely determined by choosing $l_1$, the smallest index $l$ in $L(i_1)$ such that $l\in I$ (certainly, a possible choice is to include no element of $L(i_1)$ in $I$). For instance, if $L(i_1)$ is non-empty and $l_1$ is chosen as the smallest element of  $L(i_1)$, then  $B_I=[c_{i_1},d_{l_1}]$, while if $l_1$ is not the smallest element of  $L(i_1)$, $B_I$ will be of the form $(d_{l_2},d_{l_1}]$ for $l_2$ the largest element in $L(i_1)$ smaller than $l_1$.

Let $m$ denote the cardinality of the set of indices such that $T^{dual}(x_i, t_i)$ is of the form $[c_i,\infty)$ and let $r$
be the cardinality of the set of indices such that $T^{dual}(x_i, t_i)$ is of the form $(-\infty,d_i]$. From the reasoning above, it follows that the number of choices for $B_I$ is
\[
\sum_{\mbox{choices of $i_1$}}(|L(i_1)|+1).
\]
Using that the sets $L(i_1)$ are disjoint and that the choices for $i_1$ are $m+1$ (if we count the option that all $i$ such that $T^{dual}(x_i, t_i)$ is of the form $[c_i,\infty)$ are in $I^c$), the sum above is bounded by $m+1+r$, which is bounded by $n+1$, finishing the
proof. $\square$

\subsection{Goodness of fit to multivariate normality}

In the context of testing for multivariate normality, Quiroz and Dudley, \cite{qd}, in order to establish the asymptotic distribution of their proposed procedure, considered the following class of functions on $\real^d$: Let $\H_{m}$ denote the (finite collection of) polynomials in an orthogonal basis of spherical harmonics of degree up to $m$ on the unit sphere in $\real^d$. For $h\in\H_m$, $c\in\real^d$ and $A\in$ GL($d,\real$) let
\[
g_{A,c,h}=\left\{ \begin{array}{cc}h({A(x-c)}/{\|A(x-c)\|}), & \mbox{ for } x\neq c \\
 -C, & \mbox{ for } x= c, \\ \end{array}\right.
\]
with $C$ a constant greater than $\sup_{\eta} |h(\eta)|$ (where the supremum runs over $\eta$ in the unit sphere of $\real^d$).
Let $\G=\{g_{A,c,h}: h\in\H_m, c\in\real^d, A\in \mbox{ GL}(d,\real)\}$. In \cite{qd} the metric entropy of the class $\G$ is estimated via an argument involving VC-hull classes (a concept introduced in \cite{dudc}). The uniform covering number bound obtained through this method, is the following:
\beqa\label{r10}
\mbox{For }& s=2\binom{m+d}{d}+\binom{2m+d}{d}, \mbox{ and any } v>2s/(s+1),\nonumber \\
&\sup_{Q}N_2(\eps,\G,Q)\leq K_1\exp\left(K_2/\eps^v\right).
\enqa
In particular, it was not possible to show that $\G$ was a VC-subgraph class. Now, every polynomial is
definable in the real field $(\mathbb R, +, \cdot, 0, 1)$, and so is the unitary sphere, so each of the
finite polynomials in $\H_m$ will be definable in the real field. Since multiplication, subtraction, squaring and taking square roots are definable functions, the family of the subgraphs of $g_{A,c,h}$ will be a uniformly definable family in
the real field, so by Theorem \ref{ADHMS} it will have VC-density bounded by the number of free parameters in the family, and in particular this proves that $\G$ is a VC-subgraph class.

And we can do better in computing the VC-density. Adding constants to the language does not affect o-minimality, so we can take all the parameters involved in the
polynomials in $\H_m$, add them as constants to the real field, and apply Theorem \ref{ADHMS} to this new structure. The bound
we get for the VC-density of the family of subgraphs of $g_{A,c,h}$ will be equal to the number of free parameters used in getting $A$ and $c$, so $\devc(\mbox{subgraph}(\G))\leq d^2+d$, and by Corollary \ref{Metric}, $\G$ has polynomial $L^p$ covering number with exponent $d^2+d+\eta$, for any positive $\eta$.

The large variability of the functions in $\G$ when $\|A(x-c)\|$ approaches zero, makes it difficult to apply
the method of bounded parameter space in this case. In \cite{mpq} the class $\G$ was modified, in order to avoid small values of $\|A(x-c)\|$, at the cost of sacrificing a fraction of the sample data, and only then a variation of the
bounded parameter space method was applicable. The bound given here shows that the original $\G$ is a VC-subgraph class, without need for data truncation and may help in understanding the fast convergence reported in \cite{qd} and \cite{mpq} of the finite sample distribution of the statistics proposed there to their limit distributions.

\subsection{Complexity penalties in model selection}
In Vapnik's paradigm of Structural Risk Minimization (see \cite{vv} and \cite{dgl}) in order to choose between regression models, a complexity penalty is applied to each model depending on estimates of the metric entropy of the family of functions associated. On the other hand, van de Geer \cite{vdg}, in a fairly general context, establishes the relationship between the metric entropy of classes of functions and the speed of convergence of penalized least squares estimators,  in connection with model choice. Both paradigms highlight the need for sharp estimates of metric entropy for the classes of functions defining alternative models in regression.

For example, models of the form \beq\label{r12}
Y_i=\eta(X_i^t\,\beta)+\eps_i, \qquad 1\leq i\leq n,
\enq
appear in the context of generalized linear models \cite{mn}, where $Y_i$ is the univariate response variable, $X_i$ is a $d$-dimensional vector of covariates, $\beta$ is a $d$-dimensional parameter and $\eps_i$ is the random error of the model. The function $\eta$, called the link function, is sometimes assumed to be a monotonically increasing function within a small finite set of candidates. But in a non-parametric setting (which we assume for now), $\eta$ is only required to be a continuous non-decreasing function with values in [0,1]. Thus, in the non-parametric setting, the goal is to estimate a function in $\H$, the collection of real functions on $\real^d$ of the form $\eta(x^t\,\beta)$, for $\beta\in\real^d$ and $\eta$ continuous and non-decreasing from $\real$ to $[0,1]$.

It is known that, if $\M$ denotes the collection of continuous non-decreasing functions from $\real$ to $[0,1]$, then
\[
\frac{C_{1,p}}{\eps}\leq \log \sup_QN_p(\eps,\M,Q)\leq \frac{C_{2,p}}{\eps}
\]
for $p=1,2$ and positive constants $C_{1,p}$ and $C_{2,p}$ (see the discussion in \cite{gw}). Since, clearly, covering numbers for $\H$ are larger than those for $\M$, we expect a relatively large metric entropy for $\H$, and in particular, this proves that the family of functions of the form $\eta(x^t\,\beta)$, for $\beta\in\real^d$ and $\eta$ continuous and non-decreasing from $\real$ to $[0,1]$, is not a VC-subgraph class, so one cannot expect to have any such class definable in an o-minimal structure. In fact, if one composes the increasing function $x-\sin(x)$ with any of the standard maps from $\mathbb R$ into $[0,1]$, one can easily exhibit a function which is not definable in any o-minimal expansion of the real field.

Still, in order to estimate $\eta$ (and $\beta$) non-parametrically, one could consider a sequence of nested models, as follows:
Let $\H^{(k)}$, $k\geq 2$, denote the collection of functions on $\real^d$ of the form
$\eta_k(x^t\,\beta)$,
where $\eta_k(\cdot)$ is continuous and non-decreasing from $\real$ onto [0,1] and there exist numbers $a_1<a_2<\cdots<a_k$ and $0<b_1<b_2<\cdots<b_k<1$, such that, for every $i\leq k$, $\eta_k(a_i)=b_i$, $\eta_k$ is linear between $(a_i,b_i)$ and $(a_{i+1},b_{i+1})$, for $1\leq i< k$, while for $x\leq a_1$ and $x\geq a_k$ we let
\[
\eta_k(x)=B_1 e^{c_1(x-a_1)},\>\>\mbox{ for }\>x\leq a_1\qquad\mbox{and}\qquad
\eta_k(x)=1-B_k e^{-c_k(x-a_k)},\>\>\mbox{ for }\>x\geq a_k,
\]
for positive constants $B_1,B_k,c_1$ and $c_k$, chosen to make $\eta_k$ and its derivative continuous on
the set $(-\infty,a_1]\cup[a_k,\infty)$. It seems reasonable to believe that, for moderate values of $k$, the classes $\H^{(k)}$ will provide a good approximation to an unknown function in $\H$, specially when the unknown $\eta$ is differentiable. And for the $\H^{(k)}$ the metric entropy is significantly smaller than for $\H$, as we see next. This implies (see \cite{vdg}) a much faster speed of convergence to the best approximation within each $\H^{(k)}$.

Notice first that piecewise functions in an ordered set are very easy to define if each of the components is definable. For instance, given $\bar a$, $\bar b$, $c_1, c_k, B_1, B_k$ as above, we can define the corresponding $\eta_k$ by $\eta_k(x)=y$ if and only if
\begin{footnotesize}
\[
\left(y=B_1\mbox{e}^{c_1\left(x-a_1\right)}\wedge x\leq a_1\right) \vee \bigvee_{i=1}^{k-1} \left(y=\theta(x,b_i, b_{i+1}, a_i, a_{i+1})
\wedge x\in \left(a_i, a_{i+1}\right)\right) \vee \left(y=1-B_k\mbox{e}^{c_k\left(x-a_k\right)}\wedge x\geq a_k\right),
\]
\end{footnotesize}
where $\theta(x,b_i, b_{i+1}, a_i, a_{i+1}):= \left(b_i+(b_{i+1}-b_i)\cdot \frac{x-a_i}{a_{i+1}-a_i}\right)$.

We are using $2k+2$ parameters in defining each function $\eta_k$, which means we will need $2k+d+2$ parameters to define the
function $Y_i=\eta_k(X^t\,\beta)$ (with $X$ and $\beta$ in $\mathbb R^d$). Thus, the functions $\eta_k(X^t\,\beta)$ are uniformly definable in the o-minimal model $\mathbb R_{exp}$, and, therefore, the VC-density of the family
\begin{small}
$$
\left\{
\begin{array}{rl}
  \mbox{subgraph}\left(\eta_k\left(x^t\,\beta\right)\right) : & a_1<a_2<\cdots<a_k,\> 0<b_1<b_2<\cdots<b_k<1, \nonumber \\
 &  \mbox{ with } c_1, c_k, B_1, B_k\mbox{ as defined above and }\beta\in\real^d \nonumber
\end{array}
\right\}
$$
\end{small}
is bounded by the number of parameters allowed in the definition, that is $2k+2+d$.
This implies, by our Corollary 1, that each $\H^{(k)}$ has polynomial covering number, with exponent $2k+2+d+\delta$, for any $\delta>0$.

Notice also that one can change the linear functions $\theta$, in the argument just given, for slightly more complex functions in order to guarantee any level of differentiability
at the intersections $(a_i, b_i)$ without raising the VC-density (and therefore the complexity) too much. For example, using quadratic functions instead of linear, on each interval $[a_j,a_{j+1}]$, would raise the VC-density to $3k+d+2$, and would allow us to make all the functions in $\H^{(k)}$ differentiable.

In the example we have just described and the following one it would have been somewhat unnatural to impose the assumption of total boundedness on the set of parameters defining the functions in
$\H^{(k)}$.

\subsection{Parametric estimation in Generalized Linear Models}
In the same context of generalized linear models of equation (\ref{r12}), let us move to a parametric setting by letting $\eta$ vary over all the Gaussian cummulative distribution functions, with the mean and variance, $\mu$ and $\sigma$ (as well as $\beta$), as free parameters to be estimated. We can use Corollary \ref{Metric}, to estimate the complexity of this model, as follows:
If $f$ is definable in an o-minimal expansion of $\mathbb R$, then its antiderivative (indefinite integral) belongs to the Pfaffian closure of such expansion and is therefore definable in an o-minimal structure (recall Fact \ref{Speissegger}). On the other hand,
$\exp(-x^2)$ is a definable function in $\mathbb R_{exp}$. It follows that the Gaussian density and its cumulative distribution function (c.d.f.) are both definable in $\mathbb R_{an, pfaff}$ for any choice of the parameters $\mu$ and $\sigma$, and the family of functions $\H=\{\eta(x^t\,\beta): \eta \mbox { is a Gaussian c.d.f.}, \beta\in\real^d\}$ is uniformly definable in $\mathbb R_{an, pfaff}$. Since the number of parameters involved is $d+2$, using again Corollary \ref{Metric}, we have that the class $\H$ has polynomial $L^p$ covering number with exponent $d+2+\eta$, for every $\eta>0$, uniformly on all probability laws $P$ over the pair $(X,Y)$ in (\ref{r12}).

The analysis we have just outlined would hold in exactly the same manner if, in the definition of the link function of the generalized linear model, the family of Gaussian c.d.f. is replaced by a different parametric family of distributions whose densities are uniformly definable, such as the Gamma family of distributions and others.

In future work we intend to study in more detail, the use of o-minimality methods in the context of complexity penalties for model selection.

\appendix

\section{VC-dimension vs VC-density}

In this appendix we show that, contrary to common belief (at least within
the Asymptotic Statistics community), VC-density and VC-dimension
can differ significantly over certain classes of sets.

If we restrict ourselves to finite families, it is quite easy to get
any possible difference between VC-dimension and VC-density. For
example, if we fix any $k$ points in our universe and define
$\mathcal F$ to be the family of all subsets of these fixed $k$
points, then it is easy to verify that the VC-dimension is $k$,
whereas the VC-density is 0 (the function $\Delta_\mathcal{A}(n)$ is
bounded by $2^k$ for all $n$).

\subsection{Finite unions of families of subsets}

The finite case is of course a very artificial way to force a difference
between VC-dimension and VC-density. A more common occurrence
happens when $\mathcal A$ is the union of two families. Even at the
level of the family $\mathcal F$ of semi-planes in $\mathbb R^2$, it
is easy to verify that the VC-dimension of the upper semiplanes
$\mathcal F^+$ is 2, as is the VC-dimension of the lower semiplanes
$\mathcal F^-$, whereas the VC-dimension of the union $\mathcal F$
is 3. This implies by the Sauer-Shelah Lemma\footnote{A reference to the Sauer-Shelah Lemma can be found in \cite{Sa}. An
interesting discussion
about the name is available in \cite{Bo}} that
\[
\Delta_{\mathcal F^+}(n)\leq (1/2)(n^2-n)+n+1
\]
and
\[\Delta_{\mathcal F^-}(n)\leq (1/2)(n^2-n)+n+1,\]
so by definition of $\Delta$,
\[
\Delta_{\mathcal F}(n)\leq \Delta_{\mathcal F^+}(n)+\Delta_{\mathcal
F^-}(n)\leq n^2+n+2.
\]

When taking finite unions of families, the VC-density is the maximum
of the VC-densities of the individuals in the union whereas the
VC-dimension might be increased.

This example shows a behavior that, although it happens often in the
literature, it usually never brings the difference between
VC-dimension and VC-density too far apart:

Let $\mathcal F_1, \mathcal F_2, \dots, \mathcal F_k$ each of
VC-density $N$ and suppose that $\mathcal F:=\bigcup \mathcal F_i$ has VC-density
$N+l$, so that for some set $X$ of size $N+l$ we have $|X\cap
\mathcal F|=2^{N+l}$. Now, trivially,
\[
|X\cap \F |=|X\cap \bigcup \F_i|\leq \sum_i |X\cap
{\F}_i |\leq k \sum_{j=0}^N \binom{N+1}{j},
\]
so that $k$ would need to be of the order of $2^l$. This means that
if we work with unions of $k$ families of sets, the VC-dimension
might be increased by not more than a factor of $\log(k)$.

\subsection{A bigger difference}

The final example in this section, is inspired in the finite case,
but we provide a one parameter uniformly definable family of subsets
of $\mathbb R$, with VC-dimension $N$, for any $N$.

Fix a set of $N$ points $A:=\{a_1, \dots, a_N\}$ in the interval
$(0,1)$, and for each $X_i\subset A$, let $I_i$ be a union of
subintervals of $(0,1)$ such that $I_i\cap A=X_i$ (so in particular,
we have such $I_i$ for $1\leq i\leq 2^N$). Now, let $J_i:=i+I_i$ be
the shift of the set $I_i$ by a number of units equal to its index
(so that $J_i\cap J_j=\emptyset$ for $i\neq j$), let
\[
\mathbf J:=\bigcup_{j=1}^{2^N} J_i
\]
and finally let
\[
\mathcal A:=\{ x+\mathbf{J} \}_{x\in \mathbb R}.
\]

Then $\mathcal A$ has VC-density one (by Theorem \ref{ADHMS}, since
it is a one parameter uniformly definable family in the real field)
but since for any subset $X_k\subset A$ by definition
$X_k=(-k)+\mathbf J$, we have that $|\mathcal A \cap A|=2^N$, so
that the VC-dimension of $\mathcal A$ is at least $N$, witnessed by
$A$.

\end{document}